\documentclass[12pt,a4paper]{article}

\usepackage{amsmath, amsfonts, amssymb, amsthm, calc}
\usepackage{srcltx}
\usepackage{epsfig, xspace}
\usepackage{pstricks, pst-node,pstricks-add}

\newcommand{\cF}{\mathcal{F}}
\newcommand{\cG}{\mathcal{G}}
\newcommand{\cH}{\mathcal{H}}

\newcommand{\ra}{\rightarrow}
\newcommand{\bal}{\begin{aligned}}
\newcommand{\eal}{\end{aligned}}
\newcommand{\ba}{\begin{array}}
\newcommand{\ea}{\end{array}}
\newcommand{\beq}[1]{\begin{equation}\label{#1}}
\newcommand{\eeq}{\end{equation}}
\newcommand{\N}{{\mathbb N}}

 \newcommand{\E}{{\mathbb E}}
\newcommand{\F}{{\cal F}}
\newcommand{\G}{{\cal G}}
\renewcommand{\P}{{\mathbb P}}
\newcommand{\p}{{\mathbb P}}
\newcommand{\ind}{1\hspace{-0.098cm}\mathrm{l}}

\let\BFseries\bfseries\def\bfseries{\BFseries\mathversion{bold}} 

\theoremstyle{plain}
\newtheorem{theorem}{Theorem}[section]
\newtheorem{lemma}[theorem]{Lemma}

\theoremstyle{definition}

\newtheorem{example}[theorem]{Example}

\renewenvironment{proof}[1][] {\noindent {\bf Proof#1.} }{\hspace*{\fill}$\square$\medskip\par}

\newcommand{\ie}{i.e.\@\xspace}

\author{\hspace{-2.5mm}Frank Aurzada, Hanna D\"oring, Marcel Ortgiese, and Michael {Scheutzow\thanks{Institut f\"ur Mathematik, MA 7-5, Fakult\"at II, 
        Technische Universit\"at Berlin, 
        Stra\ss e des 17.~Juni 136, 10623 Berlin, Germany;   
        {\tt \{aurzada,hdoering,ortgiese,ms\}@math.tu-berlin.de}}}}

\title{Moments of recurrence times for Markov chains}

\date{\today}

\begin{document}  \maketitle

\begin{abstract}\noindent We consider moments of the return times (or first hitting times) in a discrete time discrete space Markov chain. It is classical that the finiteness of the first moment of a return time of one state implies the finiteness of the first moment of the first return time of any other state. We extend this statement
to moments with respect to a function $f$, where $f$ satisfies a certain, best possible condition. This generalizes results of K.\ L.\ Chung (1954) who considered the functions $f(n)=n^p$ and wondered ``[...] what property of the power $n^p$ lies behind this theorem [...]''(see Chung (1967), p.\ 70). \nocite{Chung67} We exhibit that exactly the functions that do not increase exponentially -- neither globally nor locally -- fulfill the above statement.

  \par\medskip

\end{abstract}


\section{Introduction}
A classical result, see e.g.\ \cite{Kolmogorov36}, states that for any
recurrent, irreducible Markov chain on a countable state space the following holds: 
if for any state $i$ the first moment of the recurrence time is finite then this also
applies to any other state. A first generalization of this 
result appeared in~\cite{HR53}. 
If we denote by $T_{ij}$ the first time that the Markov chain visits state~$j$ 
if it is started in $i$, then the result can be stated as follows: $\E f(T_{ii}) < \infty$
for some state $i$ implies that $\E f(T_{jj}) < \infty$ for any other state $j$, 
where $f(x) = x^{n}$ for some integer $n$. The authors state the result as a lemma and
refer to the proof to an (unpublished) note by K.\ L.\ Chung and R.\ N.\ Snow. 
A  more general result can be found in~\cite[Theorem 1]{Chung54}, where $f$ is allowed to
be of the form $x^p$, for any real $p>0$. After stating the theorem, the author
also comments that the concept of generalized moments defined in terms of
a general function $f$ was suggested to him by J.\ L.\ Doob, but Chung only mentions
that his results can also be shown for functions $f$ satisfying $f(x+y) \leq f(x) + Af(y)$, for some constant~$A$.

Further related research considers recursive formulas for second moments in terms
of first moments~\cite[Sect.\ 2]{Chung54}, factorial moments~\cite{Lamperti60}, and
also, more recently, explicit formulas for higher polynomial moments~\cite{S08}.

These considerations naturally lead to the question, for which functions it is true that
a finite generalized moment of the return time for one state $i$ implies that
the moment is also finite for the return time for any other state of the Markov chain.
In this note, we characterize this class of functions.
In the following, we only consider irreducible, recurrent discrete time Markov with a countable state space.

To formulate our results, we introduce the following notation. The candidate functions $f$ 
are taken from the set
$$
\F:=\{f:\N \to (0,\infty) \mbox{ is non-decreasing and } \lim_{n \to \infty} f(n)=\infty\}.
$$
Then our objective is to classify the collection $\G$ of all $f \in \F$ such that for each irreducible recurrent discrete time Markov chain with 
a finite or countably infinite state space $E$, the following holds: 
if $\E f(T_{ii})<\infty$ for some $i \in E$ then 
$\E f(T_{jj})<\infty$ for all $j \in E$. 

Following Chung~\cite{Chung54} we additionally introduce
the class $\cH$, by stating that $f \in \cH$, if for any Markov chain the following holds:
if there exist
two states $i$ and $j$ such that $E f(T_{i j}) < \infty$ and $\E f(T_{j i}) < \infty$ then 
$\E f(T_{k\ell}) < \infty$ for any pair $k,\ell$. 

As the states $i$ and~$j$ do not have to be
distinct, it follows that $\cH$ is contained in $\cG$.
The classical result due to Kolmogorov implies that the identity function belongs
to $\cH$ and~\cite{Chung54} shows that any $f(x) = x^{p}$, for $p>0$, belongs to 
$\cH$.

Our main result states that the two classes $\cG$ and $\cH$ are in fact the same;
and we also give a characterization for a function to be in this class.

%
\begin{theorem}\label{thm:main}
Let $f \in \F$. Then the following statements are equivalent:
\begin{itemize}
 \item[(a)] $f \in \cH$;
 \item[(b)] $f \in \G$;
 \item[(c)] the following two conditions are satisfied:
	\begin{itemize}
	\item[(i)] there exists $K > 0$ such that for any $x,y> 0$, $f(x+y) \leq K f(x)f(y)$,
		
	\item[(ii)] $\limsup_{n\to\infty} \frac{1}{n} \log f(n) =0$.
	\end{itemize}
\end{itemize}
\end{theorem}

Condition (c) has an easy interpretation: It ensures that the function $f$ does not grow exponentially fast -- neither globally nor locally. In fact, one can construct functions outside $\G$ that globally increase as slowly as one wishes, but locally have parts of exponential increase (cf.\ Example~\ref{exa:constrfkt}). 

This note is structured as follows. In Section~\ref{sec:ctoa}, we prove the implication $(c) \Rightarrow (a)$ of our main theorem. In Section~\ref{sec:btoc}, we prove the implication  $(b) \Rightarrow (c)$ by showing that
any function that violates either condition $(i)$ or $(ii)$ in $(c)$ is not in the class $\G$. 
Together with the earlier
observation that $\cH \subseteq \cG$, so that $(a) \Rightarrow (b)$, this shows the equivalence
of the three statements.

\section{Proof of $(c) \Rightarrow (a)$} \label{sec:ctoa}
In this section, we prove that  $(c) \Rightarrow (a)$ in our main theorem.
We start with the following lemma, which collects some preliminary facts.

For this purpose, it is convenient to introduce the following notation: for states $i$ and $j$ of a Markov chain, we denote by $U_{ij}$ the return time from $i$ to $i$ conditioned on not crossing~$j$ (if there is such a path with positive probability). Further, $V_{ij}$ denotes the first hitting time of state $j$ when started from $i$ conditioned on not returning to $i$ before hitting $j$.

\begin{figure}[htbp]
\begin{center}
\includegraphics[scale=0.9]{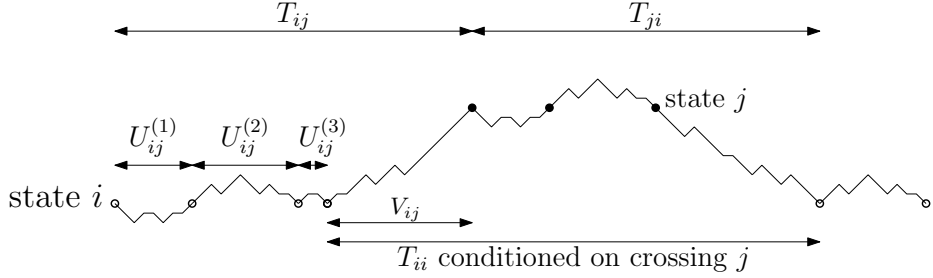}
\end{center}
\caption{Visualization of the relevant return/hitting times}
\end{figure}

\begin{lemma} \label{lem:condlm} Let $f\in \F$.
\begin{enumerate}
 \item[(i)] If for a state $i$ of a Markov chain we have $\E f(T_{ii})<\infty$ then, for any state $j$, $\E f(U_{ij})<\infty$ and $\E f(V_{ij})<\infty$.
 \item[(ii)] If for a state $i$ of a Markov chain we have $\E f(T_{ii})<\infty$ then, for any state $j$, $\E f(T_{ji})<\infty$.
\end{enumerate}
\end{lemma}

\begin{proof}
To see (i) note that if the probability $p$ of going from from $i$ to $i$ without crossing $j$ is positive, we have
$$
\infty > \E f(T_{ii}) \geq \E [ f(T_{ii}) \ind_{\{\text{do not cross $j$}\}} ] =  \E [ f(T_{ii}) \,|\, \text{do not cross $j$}\, ] \, p = p\, \E f(U_{ij}).
$$ This shows $\E f(U_{ij}) < \infty$.

Similarly, since $q:=1-p$ is the probability of first hitting to $j$ before returning to $i$, which is positive since the Markov chain is recurrent and irreducible, we have
$$
\infty > \E f(T_{ii}) \geq \E [ f(T_{ii}) \ind_{\{\text{do cross $j$}\}} ] =  \E [ f(T_{ii}) \,|\, \text{do cross $j$} \,]\, q \geq q \,\E f(V_{ij}).
$$
Here we used that $T_{ii}$ conditioned on crossing $j$ is stochastically larger than $V_{ij}$. This shows $\E f(V_{ij}) < \infty$.

To see (ii) note that
$$
\infty > \E f(T_{ii}) \geq \E [ f(T_{ii}) \ind_{\{\text{do cross $j$}\}} ] =  \E [ f(T_{ii})\, |\, \text{do cross $j$}\, ]\, q \geq q \,\E f(T_{ji}),
$$
where $q>0$ is as above. Here we used that $T_{ii}$ conditioned on crossing $j$ is stochastically larger than $T_{ji}$.
\end{proof}

We stress that, in the second part of Lemma~\ref{lem:condlm}, one cannot prove in the same way that  $\E f(T_{ij})<\infty$, since a typical path from $i$ to $i$ does not necessarily contain a path from $i$ to $j$.
However, as the next lemma shows this can be shown if we assume that condition $(c)$ holds. This lemma is also the main part of the argument for the proof of $(c) \Rightarrow (a)$ in our main theorem.

\begin{lemma} \label{lem:condlm2}  Let $f\in \F$ and assume that (c) holds.
\begin{enumerate}
 \item[(i)] If for two states $i$ and $j$ of a Markov chain we have $\E f(T_{ij})<\infty$ and $\E f(T_{ji})<\infty$ then $\E f(T_{ii})<\infty$.
 \item[(ii)] 
If for a state $i$ of a Markov chain we have $\E f(T_{ii})<\infty$ then, for any state $j$, $\E f(T_{ij})<\infty$.
\end{enumerate}
\end{lemma}

\begin{proof}
First we show (i). Clearly $T_{ii}$ is stochastically dominated by $T_{ij}+T_{ji}$, where 
$T_{ij}$ and $T_{ji}$ are independent. Using (c) and the monotonicity of $f$, we get
$$
\E f(T_{ii}) \leq \E f(T_{ij}+T_{ji}) \leq \E [ K f(T_{ij}) f(T_{ji}) ] = K \E f(T_{ij}) \E f(T_{ji}) < \infty.
$$

Now we turn our attention to (ii). For the purpose of this proof, define $f(0):=1/K$, where $K$ is as in (i) of (c). Note that $f(1)\leq f(2)\leq K f(1)^2$, so that $f(0)=K^{-1}\leq f(1)$. This shows that in fact $f(x+y)\leq K f(x) f(y)$ for {\it all} $x,y\geq 0$, since $f(x+0)=f(x)= K f(0) f(x)$.

The crucial observation (also cf.\ (\ref{eqn:lowercrucial}) below) is that
\begin{equation} \label{eqn:crucialequation}
T_{ij} = \sum_{r=1}^M U^{(r)} + V,
\end{equation}
where the random variables $U^{(r)}$ are i.i.d.\ copies of the random variable $U_{ij}$  as defined before Lemma~\ref{lem:condlm} and $V$ is a copy of the random variable $V_{ij}$  as defined before Lemma~\ref{lem:condlm}, and all variables are independent. Further, $M$ (independent of the $U$'s and $V$) is a geometric random variable with mean $1/\pi-1$, where $\pi>0$ is the probability of first hitting $j$ before $i$ when started from $i$.

It may be that $\pi=1$ -- which is the case if and only if there is no path from $i$ to $i$ without crossing $j$ -- in which case $T_{ij}=V$ and we are already done. Excluding this case, we derive from (c) and (\ref{eqn:crucialequation}) that
$$ 
\E f(T_{ij}) \leq K \E f\Big(\sum_{r=1}^{M}U^{(r)}\Big) \E f(V) = K \sum_{m=0}^\infty (1-\pi)^{m} \pi\, \E f\Big(\sum_{r=1}^m U^{(r)}\Big) \E f(V). 
$$ 
Since by Lemma~\ref{lem:condlm}(i) $\E f(V)<\infty$, it is clear that the last expression is finite provided that 
\beq{growth_sum} \limsup_{m\ra\infty } \frac{1}{m} \log \E f\Big(\sum_{r=1}^m U^{(r)}\Big) = 0\,,\eeq
where we know that $\E  f(U^{(1)})  < \infty$ from Lemma~\ref{lem:condlm}(i).

To show (\ref{growth_sum}) fix a large constant $A>0$ and estimate using (c) as follows:
\begin{align*}
 \E f(U^{(1)}+\ldots + U^{(m)}) & \leq \E f\Big( \sum_{r=1}^m (U^{(r)}\vee A)\Big)  \\
 & \leq \E f\Big( \sum_{r=1}^m (U^{(r)}\ind_{\{ U^{(r)} > A\}} + A)\Big)\\
 & \leq K f(A m) \E f\Big(\sum_{r=1}^m U^{(r)}\ind_{\{ U^{(r)} > A\}}\Big)\\
 & \leq K f(A m) \E \Big[ K^{m-1} \prod_{r=1}^m f\Big( U^{(r)}\ind_{\{ U^{(r)} > A\}}\Big) \Big]\\
 & =  f(A m) \E \Big[K f\Big( U^{(1)}\ind_{\{ U^{(1)} > A\}}\Big)  \Big]^m\, .
\end{align*}


Hence,
$$
\frac{1}{m} \log \E f(U^{(1)}+\ldots + U^{(m)})
\leq \frac{1}{m} \log f(Am) + \log (K \E f( U^{(1)}\ind_{\{ U^{(1)} > A\}})) \, .
$$
Letting $m \ra\infty$ and using part (ii) of (c) we get
$$
\limsup_{m\ra\infty } \frac{1}{m} \log \E f(U^{(1)}+\ldots + U^{(m)}) \leq \log ( K \E f( U^{(1)}\ind_{\{ U^{(1)} > A\}}))\, .
$$
Letting now $A \ra \infty$ and using dominated convergence since $\E f(U^{(1)}) < \infty$, we have
$$
\limsup_{m\ra\infty } \frac{1}{m} \log \E f(U^{(1)}+\ldots + U^{(m)}) \leq \log (K f(0)) = 0 ,
$$
which proves~(\ref{growth_sum}).
\end{proof}

Now we can prove that $(c) \Rightarrow (a)$ in our main Theorem~\ref{thm:main}.

\begin{proof}[ of $(c) \Rightarrow (a)$ in Theorem~\ref{thm:main}]
Consider a Markov chain with state space $E$, let $f$ satisfy condition~$(c)$ of Theorem~\ref{thm:main},
and assume that 
$i,j \in E$ satisfy $\E f(T_{ij}) <\infty$ and $\E f(T_{ji})<\infty$. Let $k,\ell \in E$. We want to show that 
$\E f(T_{k\ell})<\infty$.

By Lemma~\ref{lem:condlm2}(i), we have $\E f(T_{ii}) < \infty$ and $\E f(T_{jj}) < \infty$. Thus, by Lemma~\ref{lem:condlm2}(ii) $\E f(T_{j\ell})<\infty$. And by Lemma~\ref{lem:condlm}(ii) $\E f(T_{ki})<\infty$.

Therefore, using once again part (i) of (c), we have
\[
\E f(T_{k\ell}) \leq \E f(T_{ki}+T_{ij}+T_{j\ell}) \leq K^2 \E f(T_{ki}) \E f(T_{ij}) \E f(T_{j\ell})< \infty.
\]\\[-10mm]
\end{proof}

\section{Proof of $(b) \Rightarrow (c)$}\label{sec:btoc}

In this section, we will show the implication $(b) \Rightarrow (c)$ in Theorem~\ref{thm:main}
by showing that the conditions of subexponential growth rate and submultiplicativity
are in fact necessary. In Lemma~\ref{U_1+U_2}, we give abstract conditions on $f$, which
imply that $f \notin \G$, which we exploit in Lemma~\ref{lemma:sharp} to show
that if $f$ violates the submultiplicativity condition (i) of (c) in the main theorem, we have that
$f \notin \G$. Finally, in Lemma~\ref{exponential_growth}, we show that any function
growing exponentially fast, \ie does not satisfy condition~(ii), 
does not belong to $\cG$.


\begin{lemma}\label{U_1+U_2} If for a given $f \in\cF$ one can construct two independent 
random variables $U_1$ and
$U_2$ taking values in $\N$ with infinite support
such that $\E f(U_i) <\infty$ for $i=1,2$, but $\E f(U_1+U_2)  = \infty$, then
$f \notin \cG$.
\end{lemma}

\begin{proof}
Given the two random variables, we will construct a Markov chain with two special states  $0$ and $1$
with the property that $\E f(T_{11}) < \infty$, whereas $\E f(T_{00})$
is infinite, which shows that $f \notin \G$. 
The construction of the Markov chain is in some regards similar to \cite{YK39},
where the authors construct a Markov chain with $T_{00}$ having any particular distribution.

Denote by $\{x_1,x_2,\dots\}$, $\{y_1,y_2,\dots\}$ the support of $U_1$ and $U_2$, respectively. 
Formally, we can write the state space $E$ of our Markov chain as
$$
\{0\} \cup \{1\} \cup \{ (L,n,m): n\in\mathbb N, m=1,\dots, x_n-1\} \cup \{ (R,n,m): n\in\mathbb N, m=1,\dots, y_n-1\}.
$$
The state $0$ is connected only to $1$, and if in state $0$, the chain
always moves to $1$ next, \ie $p_{01} = 1$. 
If in state $1$, the chain has three possibilities. 
The first one is that the chain moves to $0$ with probability $p$, for some
parameter $p \in(0,1)$. Then, conditionally on not going to~$0$, 
with equal probability it either moves ``left'' (\ie to a state
$(L,n,1)$) or ``right'' (\ie to a state $(R,n,1)$).
\bigskip
\begin{figure}
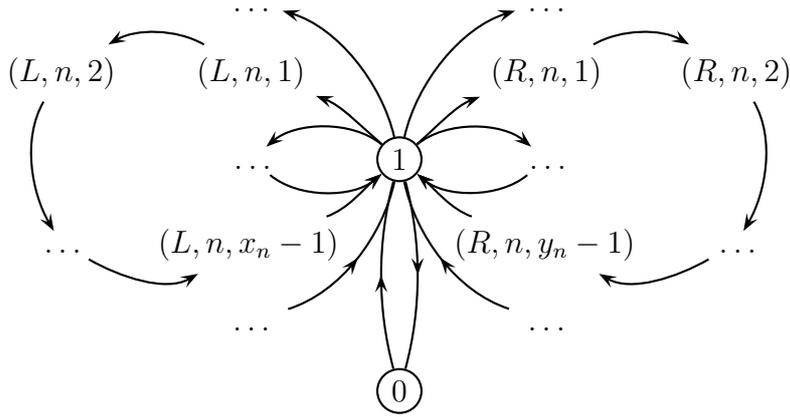

\begin{center}
$
\psmatrix[colsep=8pt, rowsep=12pt]
&{\color{white}{b}}\dots{\color{white}{b}}&&{\color{white}{b}}\dots{\color{white}{b}}&\\
{\color{white}{\Big|}}(L,n,2){\color{white}{\Big|}}
&{\color{white}{\Big|}}(L,n,1){\color{white}{\Big|}}&
&{\color{white}{\Big|}}(R,n,1){\color{white}{\Big|}}
&{\color{white}{\Big|}}(R,n,2){\color{white}{\Big|}}
&\\
&{\color{white}{g}}\dots{\color{white}{b}}
&[mnode=circle]1
&{\color{white}{g}}\dots{\color{white}{b}}&\\
{\color{white}{g}}\dots{\color{white}{b}}
&{\color{white}{\Big|}}(L,n,x_{n}-1){\color{white}{g}}
&
&{\color{white}{\Big|}}(R,n,y_{n}-1){\color{white}{b}}&{\color{white}{g}}\dots{\color{white}{b}}\\
&{\color{white}{b}}\dots{\color{white}{b}}&&{\color{white}{b}}\dots{\color{white}{b}}&\\
&&[mnode=circle]0&&\\
\psset{arrowscale=1.5}
\ncarc[arcangle=-30]{->}{2,2}{2,1}
\ncarc[arcangle=-30]{->}{2,1}{4,1}
\ncarc[arcangle=10]{<-}{2,2}{3,3}
\ncarc[arcangle=-30]{->}{4,1}{4,2}
\ncarc[arcangle=-10]{->}{4,2}{3,3}
\ncarc[ArrowInside=->, arcangle=15]{-}{6,3}{3,3}
\ncarc[ArrowInside=->, arcangle=15]{-}{3,3}{6,3}
\ncarc[arcangle=10]{->}{3,3}{2,4}
\ncarc[arcangle=30]{->}{2,4}{2,5}
\ncarc[arcangle=30]{->}{2,5}{4,5}
\ncarc[arcangle=30]{->}{4,5}{4,4}
\ncarc[arcangle=10]{->}{4,4}{3,3}
\ncarc[arcangle=-30]{->}{3,3}{1,2}
\ncarc[arcangle=30]{->}{3,3}{1,4}
\ncarc[ArrowInside=->, arcangle=30]{-}{5,4}{3,3}
\ncarc[ArrowInside=->, arcangle=-30]{-}{5,2}{3,3}
\ncarc[arcangle=-40]{->}{3,3}{3,2}
\ncarc[arcangle=-40]{->}{3,2}{3,3}
\ncarc[arcangle=40]{->}{3,3}{3,4}
\ncarc[arcangle=40]{->}{3,4}{3,3}
\endpsmatrix
$
\end{center}
\caption{Transition graph of the Markov chain}
\end{figure}

Conditionally on the next move going to the ``left'', we want $T_{11}$ to have
distribution $U_1$, therefore we set
\begin{align*}
& p_{1,(L,n,1)}= \frac{1-p}{2} \P(U_1=x_n)\, ,\\
& p_{(L,n,m-1),(L,n,m)}= 1, \quad n\in\mathbb N, 1< m \leq x_n\, ,
\end{align*}
where we identify $(1,n,x_n)$ with $1$. Similarly, conditionally on the next
move going ``right'', we would like $T_{11}$ to have the distribution
$U_2$, and set
\begin{align*}
& p_{1,(R,n,1)}= \frac{1-p}{2} \P(U_2=y_n)\, ,
\\
& p_{(R,n,m-1),(R,n,m)}= 1, \quad n\in\mathbb N, 1<m \leq y_n\, ,
\end{align*}
where we again identify $(R,n,y_n)$ with $1$.

%

Now, we can calculate the generalized moments of $T_{00}$ and $T_{11}$. 
Firstly, we find for $T_{11}$ by conditioning on the three different possibilities
\[ \E f(T_{11}) = f(2) \p \{ 1 \ra 0\} + \E f(U_1) \p \{ 1 \ra \mbox{``left''}\} + 
 \E f(U_2) \p \{ 1 \ra \mbox{``right''}\} \, ,
\]
which is finite by our assumptions on $U_1$ and $U_2$. 

However, for $T_{00}$, we obtain
\begin{equation} \label{eqn:lowercrucial}
\E f(T_{00}) = \E f\Big(\sum_{i=1}^M U^{(i)}_{10}+2\Big) \, ,
\end{equation}
where $M$ is a geometric random variable with parameter $p$ 
and $U^{(i)}_{10}$ are independent random variables that
have the same distribution
as $T_{11}$ conditioned on not going to $0$ in the first step. (Note the relation of (\ref{eqn:lowercrucial}) and (\ref{eqn:crucialequation}).)
In particular we obtain a lower bound by considering the following strategy: first 
the Markov chain jumps from $0$ to $1$, then it takes a tour to the ``left'' and after
that it takes a tour of ``right'', before it finally returns to $0$.
Thus, we obtain using that $f$ is non-decreasing
\begin{align*} \E f (T_{00}) & \geq \E[ f(2 + U_{10}^{(1)} + U_{10}^{(2)}) \, |  \mbox{ first left, then right }]\,
\p \{ \mbox{ first left, then right } \} \\
& \geq \E f(U_1 + U_2)\, p\,\big(\tfrac{1-p}{2}\big)^2 \,,
\end{align*}
where the latter is infinite by our assumptions on $U_1$ and $U_2$; and
thus as claimed $\E f(T_{00})$ is infinite.
\end{proof}

The next lemma uses the construction in Lemma~\ref{U_1+U_2}
to show that any function not satisfying the submultiplicativity condition (i) in (c) 
is not in~$\G$.

\begin{lemma}\label{lemma:sharp} Suppose that $f\in \F$ is such that for any $C > 0$, there exist
 $x_C$ and $y_C$ such that
\[ f(x_C + y_C) > C f(x_C) f(y_C) \, , \]
then $f \notin \cG$. 
\end{lemma}

\begin{proof} By Lemma~\ref{U_1+U_2} it suffices to construct two random variables $U_1$ and
 $U_2$ such that $\E f(U_i) < \infty$ for $i = 1,2$ and $\E  f(U_1+ U_2) = \infty$.

By our assumption on $f$, we can find increasing sequences $(x_k)_{k\geq 1}$ and
$(y_k)_{k\geq 1}$ such that
\[ f(x_k + y_k) > k^6 f(x_k) f(y_k) \, . \]
Indeed, to see the existence of such sequences assume that for all $x\geq x_k +1$ and $y\geq y_k+1$
the following inequality holds
\[ f(x + y) \leq (k+1)^6 f(x) f(y) \, . \]
A short calculation implies that
\[ f(x + y)
\leq (k+1)^6 \left( \frac{f(x_k\vee y_k+1)}{f(0)}\right)^2 f(x) f(y)
\quad\text{ for {\it all} }x,y\geq 0
\, , \]
in contradiction to the assumption of the lemma.

Then, for definiteness, let $U_1$ be a random variable taking value $x_k$ with 
probability $p_k:=c_1 f(x_k)^{-1} k^{-2}$ (with a suitable normalizing constant $c_1$) and similarly, $U_2$ takes
value $y_k$ with probability $q_k:=c_2 f(y_k)^{-1}k^{-2}$ (with a suitable normalizing constant $c_2$). In particular, 
we find that $\E f(U_i)$ is finite for $i=1,2$. However, 
\[ \E f(U_1+U_2) \geq \sum_{k\geq 1} f(x_k + y_k) p_k q_k \geq c_1 c_2 \sum_{k\geq 1} k^6 f(x_k) f(y_k) f(x_k)^{-1} k^{-2}f(y_k)^{-1}k^{-2} \, ,\]
which is infinite. 
\end{proof}
	
The following example exhibits a ``typical'' function that can be chosen 
to satisfy (ii) of (c), while not obeying (i) of (c).

\begin{example} \label{exa:constrfkt}
We now construct a function $f$, where we can choose the parameters in such a way that
the condition in Lemma~\ref{lemma:sharp} is satisfied and $\limsup \frac{1}{n} \log f(n) = 0$. 
We first describe the function $g = \log f$. 
Take two sequences $(s_i)_{i\geq 1}$ and $(u_i)_{i\geq 1}$ such that $u_i \leq s_{i+1} - s_i$
and $u_n \rightarrow \infty$ as $n\ra \infty$. Then, by setting $g(0) = 0$, $s_0 = u_0 = 0$
and for $i \geq 0$, define $g$ to be constant on the interval $(s_{i}+u_{i},s_{i+1})$ and for
$i \geq 1$, assume that on the interval $(s_i,s_i+u_i)$ the function $g$ grows linearly with
slope $1$.
Then, by adjusting the parameters $u_i, s_i$ (in such a way that $\lim_{n\to\infty} g(n)/n = 0$) one can make sure that the condition
in Lemma~\ref{lemma:sharp} is fulfilled for $f = e^{g}$, while at the same time by making the differences
$s_i- s_{i-1}$ large enough, one can let $f$ grow as slowly as desired.
\end{example}

The following lemma shows that the condition on the subexponential growth rate
of $f$ -- (ii) of (c) -- is really necessary. 

\begin{lemma}\label{exponential_growth} Suppose that $f\in\F$ satisfies
 \[ \limsup_{n\ra\infty} \frac{1}{n} \log f(n) > 0 \, , \]
then $f \notin \cG$. 
\end{lemma}

\begin{proof}
By our assumption on $f$ we can find an increasing sequence $x_i$ and $\delta > 0$, such that
$f(x_i) \geq e^{\delta x_i}$. Then, consider the Markov chain with two states $0$ and $1$
and transition probabilities $p_{01} = p$ and $p_{10} = 1$ for some parameter $p \in (0,1)$. 
Now, $\E f(T_{00})$ is finite, while
for any $i \geq 1$, 
\[ \E f(T_{11})  \geq e^{\delta x_i} \p \{ T_{11} \geq x_i\} = e^{\delta x_i} (1-p)^{x_i -2}\, ,
\]
which tends to infinity as $i\ra \infty$ provided that $p$ is sufficiently small.
\end{proof}

\bibliographystyle{alpha}
\bibliography{markovchains}


\end{document}